\documentclass[11pt,a4paper]{article}
\usepackage{amsmath,amssymb}

\usepackage{amsmath}

\usepackage{color}
\usepackage[colorlinks,linkcolor=blue,citecolor=red]{hyperref}

\DeclareMathSymbol{\bbbr}{\mathalpha}{AMSb}{"52}
\DeclareMathSymbol{\bbbc}{\mathalpha}{AMSb}{"52}

\newtheorem{theorem}{Theorem}

\newtheorem{definition}{Definition}

\newtheorem{lemma}[theorem]{Lemma}

\newtheorem{proposition}[theorem]{Proposition}

\begin{document}

\title{Quadratic integrals of geodesic flow, webs, and integrable billiards}

\author{{\Large Sergey I. Agafonov}\\
\\
Department of Mathematics,\\
S\~ao Paulo State University-UNESP,\\ S\~ao Jos\'e do Rio Preto, Brazil\\
e-mail: {\tt sergey.agafonov@gmail.com} }
\date{}
\maketitle
\unitlength=1mm

\vspace{1cm}

\begin{abstract}   
We present a geometric interpretation of integrability of geodesic flow by quadratic integrals in terms of the web theory and construct integrable billiards on surfaces admitting such integrals.

\bigskip

\medskip


\medskip

\end{abstract}



\section{Introduction}
This paper links two seemingly unrelated topics: the geometry of planar webs and the dynamics of geodesic flow.

 Recall that geodesic flow is a Hamiltonian system and, by the Liouville Theorem, it is enough to find only one extra first integral to describe the dynamics of geodesic flow on surfaces. If the integral is analytic then any homogeneous  part of its Taylor extension with respect to momenta  is also a first integral. This observation motivated studying homogeneous polynomial integrals since 19th century (see for instance \cite{Di-69,Da-91,Ko-96,Ko-82,Te-97,BMF-98}). 
  
The metric on a surface provides the canonical isomorphism between the tangent and cotangent spaces, which allows one to rewrite integrals, homogeneous and polynomial in momenta, also as homogeneous polynomial in velocities. Equating thus rewritten integral to zero, we get an implicit ordinary differential equation (ODE).    
It is immediate that the foliations by integral curves of this ODE are geodesic. In this paper, we study the properties of the superposition of these foliations, i.e. of the corresponding webs and nets. Note that the ODE, being implicit, may have more than one integral curve passing through a point. 

Then, for the case of a cubic integral with 3 real roots, we have a very natural interpretation of the integral: a geodesic 3-web, defined by a cubic implicit ODE, comes from a cubic integral if and only if the web is hexagonal (see \cite{Ah-19} for more detail and references concerning geodesic 3-webs).  In short, integrability means hexagonality.  

In a sense, this claim  also holds true  for quadratic integrals.  Though, one first needs  to construct 3-webs starting from a geodesic net given by integral curves of quadratic implicit ODE. The construction comes quite natural in the presence of metric structure: we complete the geodesic net  to a 4-web by the bisector net. Then, the geodesic net corresponds to a quadratic integral if and only if at least one (and therefore any) 3-subweb of the constructed web is hexagonal. Remarkably, integrability by quadratic integrals has also a conformal translation: it is equivalent to the condition that the above constructed bisector net is conformally flat (see Theorem \ref{quadrathex}).  

The configuration of a geodesic and its bisector nets brings about another interesting point. Consider  a leaf $\cal B$ of the bisector net and two geodesic rays, outbound from a point $m\in \cal B$ along the leaves of the geodesic net so that $\cal B$ does not locally separate these rays. One can view this picture either as the light ray reflecting in the mirror $\cal B$ or as the trajectory of a billiard ball hitting the wall $\cal B$ at $m$ (see the excellent introduction to the subject \cite{Ta-05}). We show (see Theorem \ref{localbilliard}) that if the geodesic net is defined by a quadratic integral then this billiard is integrable in any of the two common definitions: either as the existence of an integral of motion, globally or locally near the
the billiard wall, or as the existence of a smooth foliation  by caustics of the billiard table, globally or locally in
a neighbourhood of the billiard wall (see the survey \cite{BM-18} for the references and recent advances). The local integrability of billiards bounded by a leaf of Liouville net is not new, moreover, it was shown in \cite{IT-17} that a billiard inside a parallelogram bounded by such leaves possesses the Poncelet property (see monograph \cite{DR-11} on the Poncelet porism).  There are also globally integrable billiards with smooth boundaries on surfaces with Liouville metric (see \cite{PT-03}). We construct such billiards on surfaces with three-dimensional projective algebra and prove that the Poncelet porism holds true (see Theorem \ref{Poncelet}). 

Most of the results and considerations in this paper are local.

\section{Preliminaries} In this short section, we recall necessary definitions and facts.

\subsection{Polynomial Integrals of Geodesic Flow}
Apart from a natural geometric definition, the geodesic flow on a surface $M^2$ with a metric $g$ has  
a remarkable algebraic structure, namely, it is a Hamiltonian dynamical system on $TM^2$
with the Hamiltonian $H=\frac{1}{2}g(\xi,\xi)$. The canonical isomorphism of $T_mM^2$ and $T^*_mM^2$ via 
$g(\xi,*)=\langle f,*\rangle$, where $\xi\in T_mM^2$ is a tangent vector at a point $m$ with local coordinates $(u,v)$ and  
$f=(p,q)\in T^*_mM^2$ is the corresponding momentum, makes $p$ conjugate to $u$ and $q$ to $v$ respectively for the standard symplectic structure.

I what follows, by {\it polynomial integral of degree} $n$ we understand a first integral of geodesic flow, homogeneous in $p,q$:  
$$
I=\sum_{i=0}^na_i(u,v)p^{n-i}q^i.
$$
\begin{definition}
A  direction $[p:q]\in {\mathbb R\mathbb P^1}$ is a real root of a homogeneous polynomial $I\in \mathbb R [p,q]$, if $I(p,q)=0$.
\end{definition}
Due to the above canonical isomorphism, we can view a real root of a polynomial integral of geodesic flow as a direction in the corresponding tangent plane, i.e. we get an implicit ODE on the surface: 
$$
\sum_{i=0}^n\tilde{a}_i(u,v)du^{n-i}dv^i=0. 
$$ 
In this paper, we will be concerned with the foliations by integral curves of this ODE.    
It is immediate that these foliations are geodesic. The properties of their superposition, i.e. of the corresponding webs and nets, are studied in the next sections.  

\subsection{Hexagonal 3-webs}
Suppose that a planar 3-web ${\cal W}_3$ is given by integral curves of three ODEs 
$$
\sigma _1=0,\ \ \ \ \sigma _2=0,\ \ \ \ \sigma _3=0,\ \ \ \ 
$$
where $\sigma _i$ are differential one-forms in a planar domain. At a non-singular point, where the kernels of these forms are pairwise transverse, one can normalize the forms so that 
$\sigma_1+\sigma_2+\sigma_3=0.$ The {\it connection form} of the web ${\cal W}_3$ is a one-form $\gamma$ satisfying
$$
d\sigma_i+\gamma \wedge \sigma_i=0,\ \ \ \ i=1,2,3.
$$
The connection form depends on the normalization of the forms $\sigma_i$, however the {\it Blaschke curvature} $d\gamma$ does not.
\begin{definition}
A  3-web is hexagonal if for any non-singular point there is a local diffeomorphism mapping the web leaves in some neighbourhood of the point in 3 families of parallel line segments.   
\end{definition}
Hexagonality is equivalent to the following incidence relation, which has given its name to the notion:  for each point $m$,  any small curvilinear triangle, formed by the web leaves  and  having the point $m$ as one of its vertices, can be completed to the curvilinear hexagon, whose sides are web leaves and whose "large" diagonals also are the web leaves meeting at $m$. Computationally, both definitions of hexagonality amount to vanishing of the Blaschke curvature (see \cite{BB-38} for more detail).

\section{Dynamical proof of the Graf and Sauer Theorem }
The case of cubic integrals  has the simplest translation into the language of webs.   
\begin{theorem}\label{cubandflat}\cite{Ah-19}
A surface carries a hexagonal geodesic 3-web if and only if the geodesic flow on this surface admits a cubic first integral with 3 distinct real roots.
\end{theorem}

This theorem is an exact generalization of the Graf and Sauer Theorem \cite{GS-24} to surfaces with nonconstant gaussian curvature. To see this, observe that the flat metric   in Cartesian coordinates $x,y$ admits 3 linear integrals: $p,q$ and $qx-py$. These integrals are functionally independent. An easy exercise shows that any function of $p,q$ and $qx-py$ that is a cubic homogeneous polynomial in $p,q$, is a cubic homogeneous polynomial $I(p,q,qx-py)$ in these three integrals. The leaves of the corresponding hexagonal linear (i.e. geodesic) 3-web are integral curves of the implicit ODE $I(1,y',xy'-y)=0$, where $y'=\frac{dy}{dx}$. Applying the Legendre transform $\bar{x}=y'$, $\bar{y}'=x$, $\bar{y}+y=xy'$ to this ODE, we obtain a finite (i.e. non-differential) equation $I(1,\bar{x},\bar{y})=0$, cubic in $\bar{x},\bar{y}$.  This Legendre transform is nothing else as passing to the dual plane: $\bar{x},\bar{y}$ are the coordinates of a line, which is a solution of the ODE  $\frac{d^2y}{dx^2}=0$ for geodesics. Thus, Theorem \ref{cubandflat} provides a new, "dynamical" proof of the Graf and Sauer Theorem.

\section{Quadratic Integrals}

Consider a geodesic net $\cal G$  and define its {\it bisector net} $\cal N$ as being formed by the integral curves of directions bisecting the angles at which the geodesics of the net $\cal G$  intersect. By construction, the net $\cal N$ is orthogonal. 
\begin{definition}
An orthogonal net is called  {\it conformally flat} if it is locally diffeomorphic to one formed by coordinate lines of some isothermal coordinate system.
\end{definition}
\begin{theorem}\label{conVShex}
For the bisector net $\cal N$ of a geodesic net $\cal G$, the following conditions are equivalent:
\begin{enumerate}
\item the net $\cal N$ is conformally flat,
\item one 3-subweb of the 4-web constituted by the nets $\cal G$ and $\cal N$ is hexagonal.
\item any 3-subweb of the 4-web constituted by the nets $\cal G$ and $\cal N$ is hexagonal.
\end{enumerate}
\end{theorem}
{\it Proof:} Let us choose some local isothermal coordinates $x,y$. In these coordinates, the metric is $e^{2\lambda(x,u)}(dx^2+dy^2)$ and the equation for unparameterized geodesics assumes the form
\begin{equation}\label{eqgeoconf}
\frac{d^2y}{dx^2}=\left(\left(\frac{dy}{dx}\right)^2+1\right)\left(\lambda_y-\lambda_x \frac{dy}{dx}\right).
\end{equation}
Choose $T(x,y)$ so that the forms
$$
\omega_1=Tdx+ dy, \ \ \ \omega_2=-dx+Tdy
$$
are annihilators of the tangent vectors to the leaves of the bisector net $\cal N$. Then the annihilators of the tangent vectors to the leaves of the geodesic net $\cal G$ are
$$
u_1=\omega_1+P\omega_2, \ \ \  u_2=\omega_1-P\omega_2
$$
for a suitable function $P(x,y)$. The equation for geodesics (\ref{eqgeoconf}) gives
\begin{equation}\label{PxPy}
\begin{array}{l}
P_x=\frac{1+P^2}{P(1+T^2)}\left[ \frac{1+P^2}{1+T^2}TT_x+\frac{P^2-T^2}{1+T^2}T_y+(T^2-P^2)\lambda_x+T(1+P^2)\lambda_y \right ],\\
 \\
P_y=\frac{1+P^2}{P(1+T^2)}\left[ \frac{1-T^2P^2}{1+T^2}T_x-\frac{1+P^2}{1+T^2}TT_y+T(1+P^2)\lambda_x+(1-T^2P^2)\lambda_y \right ].
\end{array}
\end{equation}
The condition that the net $\cal N$ is conformally flat reads as (see \cite{Ac-19})
\begin{equation}\label{laplasT}
T_{xx}+T_{yy}=\frac{2T(T^2_x+T^2_y)}{1+T^2}.
\end{equation}
Finally, the Chern connections of all 3-subwebs of the 4-web formed by $\cal G$ and  $\cal N$ are equal, the exact expressions for the coefficients $\alpha,\beta $ of the connection  $\gamma=\alpha dx+\beta dy$  are given by the following formulae:
\begin{equation}
\begin{array}{l}
\alpha =\frac{1}{(1+T^2)P^2}\left[ -\frac{T^2P^2+2P^2+1}{1+T^2}TT_x+\frac{2T^2P^2+T^2+P^2}{1+T^2}T_y-T(1+P^2)(T\lambda_x+\lambda_y) \right ],\\
 \\
\beta =\frac{1}{(1+T^2)P^2}\left[ -\frac{T^2P^2+2P^2+1}{1+T^2}T_x+\frac{1-T^2P^2}{1+T^2}TT_y-(1+P^2)(T\lambda_x+\lambda_y) \right ].\\
\end{array}
\end{equation}
Now the conditions $d\gamma =0$ and (\ref{laplasT}) are equivalent modulo $d(P_xdx+P_ydy)$, where $P_x,P_y$ are as in (\ref{PxPy}), and the theorem is proved.
\hfill $\Box$
\begin{lemma}\label{realq}
If the geodesic flow on a surface admits a quadratic first integral, linearly independent of $H$ over $\mathbb R$, then it also admits a quadratic integral with real  distinct roots at a generic point.  
\end{lemma}
{\it Proof:}
Let $I_2$ be a quadratic integral, linearly independent of $H$. Then $\mu H+I_2$ is also a quadratic integral for any $\mu\in \mathbb R$. One can choose $\mu$ so that the polynomial $\mu H+ I_2$ has real roots at some neighbourhood of a generic point. If the roots are not distinct around this point then $\mu H+ I_2$ is (minus of) a square of some linear integral $L$.  Choosing  a form $\tilde{L}$,  orthogonal to $L$, we get $H=s^2 L^2 +r^2 \tilde{L}^2$ for some functions $s,r$ and obtain a quadratic integral $H-\nu L^2$ with 2 distinct real roots for any $\nu\in \mathbb R$ satisfying    $\nu > s^2$.
\hfill $\Box$
\begin{theorem}\label{quadrathex}
The geodesic flow on a surface admits a quadratic first integral if and only if the surface carries a  geodesic net $\cal G$  such that one (and therefore any) 3-subweb of the 4-web formed by $\cal G$ and its bisector net $\cal N$ is hexagonal or, equivalently, such that the net $\cal N$ is conformally flat.   
\end{theorem}
{\it Proof:}
Let $I_2$ be a quadratic integral.  We can assume that $I_2$ is chosen so that it has two real roots. Then the Liouville normal form for the metric is
\begin{equation}\label{Liouville}
g=[a^2(u)+b^2(v)](du^2+dv^2)
\end{equation}
and the quadratic integral $I_2$, written as a quadratic form, is a multiple of the "metric"
\begin{equation}\label{projective Liouville}
\bar{g}=\left[\frac{1}{a^2(u)}+\frac{1}{b^2(v)}\right]\left(\frac{du^2}{a^2(u)}-\frac{dv^2}{b^2(v)}\right),
\end{equation}
which is projectively equivalent to $g$.

The directions of the vector fields $Y_{\pm}=a(u)\partial_u\pm b(v)\partial_v$ give real roots of $I_2$. Integral curves of these vector fields form a geodesic net $\cal G$.  Since the vectors $Y_{\pm}$ have the same length, the vector fields $Y_++Y_-=2a(u)\partial_u$ and $Y_+-Y_-=2b(v)\partial_v$ are tangent to the leaves of the bisector foliations. Integral curves of these vector fields are the coordinate lines of the isothermal coordinates $(u,v)$. By Theorem \ref{conVShex}, this condition is equivalent to hexagonality of any 3-subweb formed by $\cal G$ and $\cal N$.

To prove the converse claim, choose isothermal coordinates $u,v$ such that the bisector net $\cal N$ of the given geodesic net $\cal G$ is formed by the coordinate lines. Then the foliations of $\cal G$ are integral curves of the ODEs 
$$dv-P(u,v)du=0\ \ \ \ {\rm and} \ \ \ \  dv+P(u,v)du=0$$
 for some function $P(u,v)$. The equation for unparameterized geodesics takes the form
\begin{equation}\label{eqgeoconf2}
\frac{d^2v}{du^2}=\left(\left(\frac{dv}{du}\right)^2+1\right)\left(\frac{\Lambda_v}{2\Lambda}-\frac{\Lambda_u}{2\Lambda} \frac{dv}{du}\right),
\end{equation}
where $\Lambda(u,v)$ is fixed by the metric $g=\Lambda(u,v)(du^2+dv^2)$. Equation (\ref{eqgeoconf2}) gives 
$$
P_u= - \frac{P(1+P^2)}{2}\frac{\Lambda_u}{\Lambda}, \ \ \ \ \ \ \ P_v=  \frac{(1+P^2)}{2P}\frac{\Lambda_v}{\Lambda}. 
$$ 
Now the condition $d(P_udu+P_vdv)=0$ yields $\Lambda_{uv}=0$. Therefore the metric has Liouville's form and its geodesic flow admits a quadratic integral.
\hfill $\Box$\\

\medskip
\noindent {\bf Remark 1.} It is not difficult to find an explicit formula for the quadratic first integral $I_2$ whose real roots coincide with the slopes of the tangent spaces to the leaves of $\cal G$. This integral has the form $I_2=\mu(u,v) (q-Pp)(q+Pp)$. Splitting the equation $\{H,I_2\}=0$ by monomials $p^jq^{3-j}$, $j=0,...,3$, we get a system of PDEs for the integrating factor $\mu$. This system can be easily integrated to give $\mu=\frac{1}{1+P^2}.$\\

\noindent {\bf Remark 2.} It is well known (see \cite{BB-38}) that a 4-web, whose all 3-subwebs are hexagonal, is linearizable and the foliations of its "linear form" are 4 pencils of lines. The  4-web constituted by the nets $\cal G$ and $\cal N$ is even more degenerate. Namely, the 4 pencils of its linearization are collinear and the pencil vertices form harmonic quadruplet. In fact, the first integrals of the foliations corresponding to the zero directions of (\ref{projective Liouville}) can be chosen as $I_{\pm}=\int \frac{du}{a(u)}\pm \int \frac{dv}{b(v)}$. The foliations of the bisector net $\cal N$ have first integrals $u$ and $v$. Thus, in the coordinates $U=\int \frac{du}{a(u)}$ and  $V=\int \frac{dv}{b(v)}$, the foliations of our 4-web are the lines $U=const$, $V=const$, $U\pm V=const$. Locally, the picture is exactly the same as for the 4-web formed by confocal conics and tangent lines to one conic of the confocal family (see \cite{Ac-19}).  \\

\section{Linear Integrals}
For completeness, let us consider the foliation constructed by integration of the zero direction of a linear integral.
By Noether's Theorem, the existence of a linear integral of the geodesic flow on a surface is equivalent to the existence of a Killing vector field $X$, which is the projection of the Hamiltonian vector field corresponding to the integral.
Suppose that the Killing vector field does not vanish at $m_0\in M$. Then  we can "rectify" it in some neighborhood of $m_0$ by choosing local coordinates $u,v$: $X=\partial_u$. In these coordinates, the coefficients of the metric do not depend on $u$. Applying a suitable substitution $v\to \psi(v)$, $u\to u+\varphi(v)$
we bring the metric to the form $g=E(v)du^2+dv^2$.
Now, the integral is $I_1=p$ and the ODE for its zero direction reads as $du=0$. Thus, we have the following geometric observation. \begin{proposition}\label{linearKilling}
The geodesic flow on a surface admits a linear first integral if and only if the surface carries a  geodesic foliation  such that its orthogonal foliation is formed by orbits of a Killing vector field.
\end{proposition}

\section{Quadratic Integrals and Integrable Billiards} 
Utilizing the obtained geometric  interpretation of quadratic integrals, one can construct locally integrable billiards, mimicking the classical billiards in ellipses in the Euclidean plane. 
Integrability here is present in any of the two common definitions: either as the existence of an integral of motion, globally or locally near the
billiard wall, or as the existence of a smooth foliation  by caustics of the billiard table (or locally in
a neighbourhood of the billiard wall). 

The construction is as follows. 

\begin{theorem}\label{localbilliard}
Suppose there is  a quadratic integral of the geodesic flow on a surface such that 
\begin{enumerate}
\item its zero directions are real,
\item the geodesic net defined by these directions has a closed simple smooth envelope,
\item the foliation of the bisector net that goes "along" the envelope has a closed simple smooth leaf $\mathfrak{B}$ convex with respect to the geodesics,  
\end{enumerate}
then the billiard with the wall $\mathfrak{B}$ is locally integrable, namely 
\begin{enumerate}
\item there is an integral of motion for trajectories forming a  sufficiently small angle with the billiard wall, and
\item a neighbourhood of the billiard wall is foliated by caustics.
\end{enumerate}
\end{theorem}     

{\it Proof:} In some local coordinates, the metric has the normal Liouville form (\ref{Liouville}). Then it is well known that the geodesic flow has the following first integral related to the projectively equivalent metric (\ref{projective Liouville}):
$$
I_0=\left(\frac{\det(g}{\det(\bar{g})}\right)^{\frac{2}{3}}\bar{g}=[a^2(u)+b^2(v)](b^2(v)du^2-a^2(u)dv^2). 
$$

This integral  has real zero directions. Without loss of generality, we can assume that $I_0$ is the chosen integral. The envelopes of the corresponding geodesic net are the coordinate lines $u=u_0$ and $v=v_0$, where $u_0$, $v_0$ are zeros of $a$ and $b$ respectively. Suppose that the envelope satisfying the hypothesis is $v=v_0$. Then the foliation $\mathcal{F}$ of the bisector net, whose leaves goes along the chosen envelope, are the one by coordinate curves $v=const$. Let the leaf $\mathfrak{B}$ satisfying the hypothesis be $v=v_1$.  

       The integral $I_0$ generates  the pencil of quadratic integrals 
$$
I_{\mu}=\mu g+I_0=[a^2(u)+b^2(v)]([b^2(v)+\mu]du^2-[a^2(u)-\mu]dv^2),
$$ each integral of the pencil having real zero directions for sufficiently small $\mu$.  
Thus, we have a one-parameter family of geodesic nets $\cal G_{\mu}$, whose leaves are the integral curves of the ODE 
$$
\frac{du}{\sqrt{a^2(u)-\mu}} \pm \frac{dv}{\sqrt{b^2(v)+\mu}}=0.
$$ 
The envelopes of the corresponding  geodesic foliations are again the coordinate lines $v=const$ (and $u=const$). 
Invoking the Liouville formula for geodesics, we conclude that any geodesic is tangent to some envelope   
labeled by $\mu$. Taking into account that the wall $\mathfrak{B}$ is a leaf for the bisector net for each geodesic net $\cal G_{\mu}$, we conclude that $\mu$ is the integral of motion and that  the envelopes $v=const$ are the caustics.   
\hfill $\Box$\\

Apart from the billiards in ellipses, there are other examples of {\it globally} integrable billiards  resulting from the above described construction (cf. \cite{PT-03}). Here we describe examples  on surfaces with three-dimensional algebra of projective symmetry. (Recall that an infinitesimal projective symmetry is a vector field whose local flow respects the unparametrized geodesics.) We will follow the description of geodesics given in \cite{Ah-19}. Our surface is the upper half-plane $M^2=\{(z,y)\in \mathbb{R}^2:z>0\}$, the equation for unparametrized geodesics is 
\begin{equation}\label{eqgeo3}
z^3\frac{d^2y}{dz^2}=\epsilon\left(\frac{dy}{dz}\right)^3,\ \ \ \epsilon=\pm 1.
\end{equation}
There is a parametric family of metrics sharing these geodesics (see \cite{BMM-08}). 
The equation for geodesic can be easily integrated to give the general solution in the implicit form
\begin{equation}\label{solgeo3}
k^2(y-l)^2-kz^2=\epsilon,
\end{equation}
where $k,l$ are integration constants. 
Obtained by limit $k\to \infty$, there are special solutions $y=const$,  which will play an important role. 
Thus, we have a 2-parametric family of conics of the form
\begin{equation}\label{incidence}
Ay^2+2By+C+Dz^2=0,
\end{equation}
the special solutions being considered as double lines $(y-l)^2=0$.
The closure of this family, called the {\it dual} space of $M^2$ in what follows, is the quadric 
\begin{equation}\label{quadric}
AC-B^2+\epsilon D^2=0
\end{equation}
in $\mathbb R\mathbb P^3$.
The special solutions correspond to the section of the above quadric by the plane $D=0$. This section is a smooth conic $c_0$. To a point $(z,y)\in M^2$ on the surface, there corresponds the section of the quadric (\ref{quadric}) by the plane (\ref{incidence}), intersecting the plane $D=0$ along a line tangent to the conic $c_0$. 

To compute quadratic integrals, one needs a metric, but the corresponding geodesic nets are the same for all metrics with geodesics (\ref{solgeo3}). These nets are described as follows. 
\begin{theorem}\label{geonets}
Geodesics of the net, determined by a quadratic integral of the geodesic flow of any metric with geodesics (\ref{solgeo3}), are points on some section of the quadric (\ref{quadric}) by a plane different from the plane $D=0$.     
\end{theorem} 
{\it Proof:} The proof is based on the description of hexagonal geodesic 3-webs on our surfaces given in \cite{Ah-19}: the geodesics of one foliation of the web trace an arc of the conic  $c_0$ and the geodesics of the other two run along arcs of a section of the quadric (\ref{quadric}) by some fixed plane, different from the plane $D=0$. The cubic integrals of geodesic flow on our surface are products of a linear and of a quadratic integral (see, for instance, \cite{MS-11}). The geodesic foliation, corresponding to $c_0$, is the same for all hexagonal geodesic 3-webs. Therefore it is the one of the linear integral. Hence the theorem claim. 
\hfill $\Box$\\

Let us now represent the construction of billiard via quadratic integral in the dual space. 
The integral fixes a conic $c_I$ on the quadric (\ref{quadric}). The billiard trajectory is the sequence of segments of geodesics $l_i\in c_I$ joined by conics $c_i$ on the quadric (\ref{quadric}). The conics $c_i$  are tangent to the conic $c_0$. 
\begin{lemma}
For any billiard trajectory with the caustic determined by $c_I$, the conics $c_i$ are tangent to some 
plane conic $c_w\ne c_0$ on the quadric (\ref{quadric}).   
\end{lemma}    
{\it Proof:} To the conics $c_i$ on the quadric (\ref{quadric}), there correspond points   $m_i\in M^2$. The points lie on the billiard wall $\mathfrak{B}$. Recall that the wall is also an envelope of the foliation leaves of some geodesic net $\cal G$ coming from a quadratic integral. Therefore all the points $m_i$ are incident to geodesics of $\cal G$, forming a conic $c_w$ in the dual space. Thus, as $m$ ranges over $\mathfrak{B}$, the corresponding conics in the dual space envelope $c_w$.\hfill $\Box$\\

The geometric construction in the proof of the following theorem  implies that our billiard is globally integrable. 
\begin{theorem}\label{Poncelet}
Poncelet's porism holds true for any billiard constructed via a quadratic integral on a surface with three-dimensional algebra of projective symmetries.  
\end{theorem}
{\it Proof:} We reduce the theorem claim to the one for billiards in ellipses in the Euclidean plane.
Consider a point $m$ on the wall $\mathfrak{B}$. The corresponding conic $c_m$ on the quadric (\ref{quadric}) is tangent to $c_0$ and to $c_w$ at points $q_0(m)$ and $q_w(m)$. As $m$ runs along $\mathfrak{B}$ the lines $q_0(m)q_w(m)$ spread a ruled surface $K\subset \mathbb R\mathbb P^3$. This surfaces is the envelope of one of two possible smooth families of planes rolling on $c_0$ and $c_w$. 
\begin{lemma}\label{alglemma}
The ruled surface $K$ is a cone. 
\end{lemma}
Now the image of a billiard trajectories in the dual space is the following. Each conic $c_i$ of the sequence of conics $\{c_i\}_{i\in \mathbb N}$ intersects the conic $c_I$, corresponding to the integral, at points $l_i$, corresponding to geodesic segments of the trajectory,  so that $l_i \in c_i\cap c_{i+1}$. If we project the quadric from the vertex $v$ of the cone $K$ into some plane then we get the dual image of a trajectory for a billiard inside a conic in the Euclidean plane, and  the Poncelet's  porism follows. \hfill $\Box$

\medskip   
 {\it Proof of Lemma \ref{alglemma}:} The planes of $c_0$ and $c_w$ intersect along a line $L_w\in\mathbb R\mathbb P^3$. Consider a point $q\in L_w$. The tangent lines from $q$ to the quadric (\ref{quadric}) touch the quadric along a plane conic $c_q$, which is the intersection of the polar plane $\pi_ q$ of $q$ with the quadric.  There are two tangent lines from $q$ to the quadric that are tangent to $c_0$ and two lines that are tangent to $c_w$. Two of these 4 lines are $qq_0(m)$ and $qq_w(m)$ for some $m\in \mathfrak{B}$. Let the other pair be $qq_0(\bar{m})$ and $qq_w(\bar{m})$, where the point $\bar{m}$ is some other point on the wall $\mathfrak{B}$. Choose the point $q\in L_w$ so that the above described  4 points $q_0(m)$, $q_w(m)$, $q_0(\bar{m})$, and $q_w(\bar{m})$ are distinct.
 The lines $q_0(m)q_w(m)$ and $q_0(\bar{m})q_w(\bar{m})$, being coplanar, intersect at some point $v\in\mathbb R\mathbb P^3$. Consider the cone over $c_0$ with the vertex $v$. This cone intersect the quadric (\ref{quadric}) along a curve of degree 4. This curve decomposes into the conic $c_0$, which is of degree 2, and some curve $c$ that must be also of degree 2. Hence $c$ is a conic, containing the points $q_w(m)$ and $q_w(\bar{m})$ and tangent to the lines $qq_0(\bar{m})$ and $qq_w(\bar{m})$. Since the conic $c$ lies on the quadric (\ref{quadric}), it is a plane section of the quadric. But the section plane is fixed by the lines $qq_0(\bar{m})$ and $qq_w(\bar{m})$. Thus $c=c_w$ and hence the lemma claim. 
\hfill $\Box$\\

\section{Concluding remarks} The "hexagonality property" of integrable billiards, established in this paper, naturally supplements four equivalent properties of integrable billiards \cite{GIT-19} by a fifths one.

\section*{Acknowledgements}
This research was supported by FAPESP grant \#2018/20009-6 and partially by DAAD grant \#91727921.
The author thanks V.S. Matveev and I. Izmestiev for useful discussions.


\begin{thebibliography}{99}

\bibitem[Ac-19]{Ac-19} S.I. Agafonov, Confocal conics and 4-webs of maximal rank, arXiv: 1912.01817 [math.DG]

\bibitem[Ah-19]{Ah-19} S.I. Agafonov, Hexagonal geodesic 3-webs, {\it Int. Math. Res. Not. IMRN,} (2019), 33 p. \\
DOI: 10.1093/imrn/rnz172   arXiv: 1806.03072 [math.DG]

\bibitem[BM-18]{BM-18} M. Bialy, A.E. Mironov, 
A survey on polynomial in momenta integrals for billiard problems.  
{\it Philos. Trans. Roy. Soc.} A 376 (2018), no. 2131, 19 pp. 


\bibitem[BB-38]{BB-38} W. Blaschke,  G. Bol,
{\it Geometrie der Gewebe}, Topologische Fragen der
Differentialgeometrie. J. Springer, Berlin, 1938.

\bibitem[BMF-98]{BMF-98} A.V. Bolsinov, V.S. Matveev, A.T. Fomenko, Two-dimensional Riemannian metrics with an integrable geodesic flow. Local and global geometries, {\it Mat.
Sb.} 189 (10) (1998) 5--32. Engl. translation: {\it Sb. Math.} 189 (9–10) (1998) 1441--1466.

\bibitem[BMM-08]{BMM-08} R.L. Bryant, G. Manno, V.S. Matveev,
A solution of a problem of Sophus Lie: normal forms of two-dimensional metrics admitting two projective vector fields.
{\it Math. Ann.} 340 (2008), no. 2, 437-463.


\bibitem[Da-91]{Da-91} G. Darboux, { \it Le\c cons sur la th\'eorie g\'en\'erale des surfaces et les applications g\'eom\'etriques du calcul infinit\'esimal.} Paris: Gautier, Villar, 1891.

\bibitem[Di-69]{Di-69} U. Dini, Sopra un problema che si presenta nella teoria generale delle rappresentazione geografiche di una superficie su di un’altra. (Italian)
{\it Brioschi Ann.} (2) III,  (1869) 269--294.


\bibitem[DR-11]{DR-11} V. Dragovi\'c, V. Radnovi\'c, {\it Poncelet porisms and beyond. Integrable billiards, hyperelliptic Jacobians and pencils of quadrics.} Frontiers in Mathematics. Birkh\"auser/Springer Basel AG, Basel, 2011.

\bibitem[GIT-19]{GIT-19}A. Glutsyuk, I. Izmestiev, S. Tabachnikov, Four equivalent properties of integrable billiards.
arXiv:1909.09028 [math.DS]

\bibitem[GS-24]{GS-24} H. Graf, R. Sauer,  \"Uber dreifache
Geradensysteme in der Ebene, welche Dreiecksnetze bilden, {\it
Sitzungsb. Math.-Naturw. Abt.} (1924), 119-156.

\bibitem[IT-17]{IT-17} I. Izmestiev, S. Tabachnikov,
Ivory's theorem revisited.
{\it J. Integrable Syst.} 2 (2017) no. 1, 36 pp.

\bibitem[Ko-96]{Ko-96} G. Koenigs, Sur les gé\'eodesiques a int\'egrales quadratiques, in: Note II from Darboux {\it Leçons sur la th\'eorie g\'en\'erale des surfaces}, vol. IV, Chelsea
Publishing, 1896.

\bibitem[Ko-82]{Ko-82} V.N. Kolokol'tsov, Geodesic flows on two-dimensional manifolds with an additional first integral that is polynomial with respect to velocities, {\it Izv.
Akad. Nauk SSSR Ser. Mat.} 46 (5) (1982) 994--1010, 1135. Engl. Transl.: {\it Math. USSR-Izv} 21 (1983) 291--306.

\bibitem[MS-11]{MS-11} V.S. Matveev, V.V. Shevchishin, Two-dimensional superintegrable metrics with one linear and one cubic integrals, {\it ‎J. Geom. Phys.} 61 (2011) 1353--1377.

\bibitem[PT-03]{PT-03} G. Popov, P. Topalov, Liouville billiard tables and an inverse spectral result. {\it Ergodic Theory Dynam. Systems} 23 (2003), no. 1, 225--248. 

\bibitem[Ta-05]{Ta-05} S. Tabachnikov, {\it Geometry and billiards.}  American Mathematical Society, Providence, RI, 2005.

\bibitem[Te-97]{Te-97} V.V. Ten,
Local integrals of geodesic flows. (Russian)
{\it Regul. Khaoticheskaya Din.} 2 (1997) no. 2, 87--89.



\end{thebibliography}
\end{document}